# ON AN APPLICATION OF THE SUSLIN MONIC POLYNOMIAL THEOREM (I) .

**C.L.Wangneo**

**Jammu,J&K,India,180002**

**(E-mail:-wangneo.chaman@gmail.com )**

**Abstract:-**

In this paper we state results that lead to our main theorem which states the following ;

Main Theorem :- Let $B = D[x_1,-,x_n]$, be a polynomial ring in the commuting variables $x_i$ over a division ring D. Let M be a finitely generated B-module. Let $B'_m = D[x_1,-,x_m]$, be the polynomial ring in 'm' variables; where m is an integer such that $0 \leq m \leq n$, with $B'_0 = D$, and $B'_n = B$. Then the following conditions on M are equivalent ;

(i) Krull dimension $(M) = m$; $0 \leq m \leq n$ .

(ii) M is a non -torsion $B'_m$ -module (not necessarily finitely generated as a $B'_m$ -module ) and such that for any k>m , M is a torsion $B'_k$ - module .

We then also announce a generalisation of the above theorem .

**Introduction**

In this paper we state results that culminate in our main theorem which is a characterisation of a finitely generated module M over a polynomial ring $B = D[x_1,-,x_n]$, over a division ring D in the commuting variables $x_i$ with krull dimension (M)=m . The statement of our main theorem is the following;

Main Theorem :- Let $B = D[x_1,-,x_n]$, be a polynomial ring in the commuting variables $x_i$ over a division ring D. Let M be a finitely generated B-module. Let $B'_m = D[x_1,-,x_m]$, be the polynomial ring in 'm' variables; where m is an integer such that $0 \leq m \leq n$, with $B'_0 = D$, and $B'_n = B$. Then the following conditions on M are equivalent ;

(i) Krull dimension (M) = m; $0 \leq m \leq n$ .

(ii) M is a non -torsion $B'_m$ -module (not necessarily finitely generated as a $B'_m$ -module ) and such that for any k>m , M is a torsion $B'_k$ - module .

We then also announce a generalisation of the above theorem .

The paper is divided into two sections . In section (1) we first state the Suslin monic polynomial theorem which is verifiable on the same lines as given in [3] and which holds for the polynomial ring $B = D[x_1,-,x_n]$, over a division ring D in the commuting variables $x_i$ . We then introduce definitions and results in section (1) that yield theorem- 0 for the polynomial ring $B = D[x_1,-,x_n]$ . Our theorem-0 states the following ;

**Theorem-0**  :- Let  B be the Polynomial ring $B= D[x_1,-,x_n]$,  where D is a division ring  and $x_i$ are commuting variables over D. Let M be a finitely generated B-module. Then the following  hold true for  M  ;

(a) If M is a torsion B-module, then there exists a change of variables $t_i$ of $x_i$ such that B can also be written as a polynomial ring $B = D [t_1,-, t_n]$, in the new variables $t_i$ and such that M is a finitely generated $B_{n-1}$-module, where $B_{n-1}$ is the polynomial ring  $B_{n-1} = D[t_1,\ldots,t_{n-1}]$  in the fewer new   'n-1' variables.

(b) In case M is further a torsion $B_{n-1}$-module ($B_{n-1}$ is as in (a) above) then there exists a further change of variables $z_i$ of $t_i$; $1 \le i \le n-1$, such that $B= D[z_1,\ldots z_{n-1}]$  and M is a finitely generated  $B_{n-2}$-module where  $B_{n-2}$ is the polynomial ring  $B_{n-2} = D [z_1,-,z_{n-2}]$. Moreover in this case we can write B also as the polynomial ring   $B = D[z_1,-,z_n]$,   where $z_n$ can be chosen such that   $z_n=t_n$ .

An application of   theorem-0  makes it possible to state theorem -I  in section (1)  and  which  makes the statement of  theorem-0    more precise by connecting it with the Krull dimension of the module M . Our theorem-I states  the  following ;

Theorem-I :- Let  M  be  a  finitely  generated B-module where  B is the Polynomial ring   $B = D[x_1,-,x_n]$,    in n-commuting variables $x_i$, over a division  ring D . Let   $B'_m$  $=D[x_1,-,x_m]$, be  the polynomial ring in 'm' variables; where m is an integer such that $0 \le m \le n$, with  $B'_0$ =D, and

$B'_n$ =B. Then the following  conditions  on  M  are  equivalent ;

(i) M  is  a  finitely generated B-module ,  with $|M|_B$ =m ; $0 \le m \le n$ .

(ii)There exists a change of variables $t_i$ of $x_i$ such that B can also be expressed as a Polynomial ring $B = A[t_1,-,t_n]$, in the new variables $t_i$ ( giving an automorphism f of B such that $f(x_i)=t_i$ ) and M is a finitely generated non-torsion $B_m$ -module where $B_m$ is the polynomial ring $B_m = D[t_1,-,t_m]$, in the fewer, m variables , $t_1,-,t_m$ ( clearly $f(B'_m)=B_m$ ) . (Moreover , in this case then for any positive integer k , with k>m and $0 \le k \le n$ , M is a finitely generated torsion $B_k$ module where $B_k$ is the polynomial ring $B_k = D[t_1,-,t_k]$, in the fewer, k variables , $t_1,-,t_k$ ) .

 A straight-forward version of the above theorem then yields our main theorem as stated in the beginning .

Then in section (2) we introduce further definitions and results which make it possible for us to state and announce the generalisations of the key results of section (1) .

**Notation and Terminology:-**

Throughout in this paper a ring is meant to be an associative ring with identity which is not necessarily a commutative ring. We will throughout adhere to the same notation and definitions as in [2] . Thus by a noetherian ring R we mean that R is both a left as well as a right noetherian ring. By a module M over a ring R we mean that M is a right R-module unless stated otherwise. For the basic definitions regarding noetherian modules and noetherian rings and all those regarding Krull dimension, we refer the reader to [2]. Also if M is a right noetherian module over a right noetherian ring B , then we will denote by $|M|_B$ the right krull dimension of M . $|M|$ also denotes the right krull dimension of M in case there is no confusion regarding the ring over which we are considering M as a right module . Hence if R is a ring with krull dimension then we denote the krull dimension of R , usually , by $|R|$ . If R is a ring, then the polynomial ring R[x] over R is the usual polynomial ring over R in a commuting variable x. Throughout, in this paper, by a polynomial ring, say $R[x_1,-,x_n]$, is meant to be a polynomial ring in several variables $x_1,-,x_n$ over R that commute with each other and that also commute with the elements of R . Lastly, we mention that since we treat the key results of section (1) as natural consequences of the statements of the intermediary results , thus we do not include any proofs in this paper .

**Section(1)  (Main Theorem)**  :- As stated in the introduction , in this section we state definitions and results that lead to our main theorem for the polynomial ring B= $D[x_1,-,x_n]$, over a division ring D in the 'n' commuting variables $x_i$ . For this we first state the Suslin monic polynomial theorem , which can be verified on the same lines for the polynomial ring B= $D[x_1,-,x_n]$, over a division ring D in the 'n' commuting variables $x_i$ , as given in [3] . We state this theorem below;

**Theorem (1.1) (Suslin Monic Polynomial Theorem):-**

Let B be a polynomial ring B = $D[x_1,-,x_n]$, over a division ring D in the 'n' commuting variables $x_i$. Then the following holds true ;

If J is a nonzero right ideal of B such that J contains a regular element of B, then there exists a change of variables, say $t_i$ of $x_i$, such that B can also be written as a polynomial ring B = $D[t_1,-,t_n]$, in the variables $t_i$ and such that J contains a polynomial f written in terms of these new variables $t_i$, which is monic in $t_n$. Moreover in this case B/J is a finitely generated $B_{n-1}$-module where $B_{n-1}$ is the polynomial ring $B_{n-1} = D[t_1,-,t_{n-1}]$, in the 'n-1' variables $t_i$.

We next modify the above theorem to yield theorem-0 . For this we further state two lemmas that are needed for the necessary modification as below ;

**Lemma (1.2)** :- Let $B = D[x_1,-,x_n]$, be a polynomial ring in n-commuting variables $x_i$ over a division ring D. Let $t_i$ be a change of variables of $x_i$ such that for an integer $d \geq 1$.

$$t_i = \{ x_n - x_1{}^d, \quad \text{if } i=1 \} \}$$
$$= \{ x_i, \text{if } 1 < i < n \} \qquad (1)$$
$$= \{ x_1, \text{if } i=n \quad \}$$

Then we can write , $B = D[t_1,-,t_n]$, a polynomial ring in the variables $t_i$ . Converse is also true, namely if $B = D[t_1,-,t_n]$, is a polynomial ring in the commuting variables $t_i$ over D and $t_i$ are as in (1) above , then $B = D[x_1,-,x_n]$, is a polynomial ring in the variables $x_i$.

**Lemma (1.3):-** Let everything be as in Lemma (1.2) above and let

$B = D[x_1,-,x_n]$, be a polynomial ring in the commuting variables $x_i$ over D. Then $B = D[t_1,-,t_n]$, where $t_i$ is a change of variables $x_i$ given by (1) of Lemma (1.2) above . Moreover if the ring $B_{n-1}$ denotes the ring $B_{n-1} = D[t_1,-,t_{n-1}]$, and if $z_i$ is a change of variables $t_i$, $1 \leq i \leq n-1$ such that for some integer $l \geq 1$;

$z_1, = t_{n-1} - t_1{}^l$ , $z_{2=}t_{2, .........}$ , $z_{n-2} = t_{n-2}$ and $z_{n-1} = t_1$

 Then we have the following:

(1) $B_{n-1} = D[z_1,-,z_{n-1}]$, is a polynomial ring in the 'n-1' commuting variables $z_i$ .

(2)If $z_n = t_n = x_1$, then $z_1, z_2, \ldots, z_{n-1}$, $z_n$ are also algebraically independent over D and we have $B = D[z_1, \text{-}, z_{n-1}, z_n] = D[t_1, \text{-}, t_n] = D[x_1, \text{-}, x_n]$.

**Remark:-**Note that in the above two lemmas we have been able to produce all in all n- variables $z_1, \ldots, z_n$ which is a common change of variables $x_i$ and $t_i$ over B such that $B = D[z_1, \text{-}, z_n]$ . Moreover if $B_{n-1} = D[z_1, \text{-}, z_{n-1}]$, and $B_{n-2} = D[z_1, \text{-}, z_{n-2}]$ , then we can apply the above stated lemmas to the study of a module M which is a finitely generated, torsion B - module as well as a torsion $B_{n-1}$-module , respectively as given in our theorem-0 stated below.

However before we state our theorem-0 we first introduce some well known definitions for the benefit of the reader .

**Definition (1.4):-** Given a noetherian ring R, denote by C(0) the set of regular elements of R. Let M be a module over R. An element, y of M is said to be a torsion element of M if $yd = 0$ for some element d of C (0). y is said to be a torsion-free element of M if $yd \neq 0$ for all elements d of C(0). Clearly y is a torsion-free element implies that y is a nonzero element of M.

**Definition (1.5):-** We say a module M over a noetherian ring R is a torsion module if every element of M is a torsion element. We say a module M over a noetherian ring R is a torsion-free module if every nonzero element of M is a torsion-free element.

**Remark:** - Let R be a right noetherian ring and let M be a finitely generated module over R. Let T(M) be the set of torsion elements of M. Then T(M) is a non-empty subset of M and T(M) is usually not a submodule of M. However T(M) is a submodule of M if C(0), the set of regular elements of R is a right ore set. If R is a semiprime right Noetherian ring then by Goldie's Theorem (see [1]) C(0) is a right ore subset of R.

Now we mention that a judicious use of the Suslin monic polynomial theorem ,namely theorem (1.1) , stated above yields our theorem (1.6) ( named as theorem-0) for the polynomial ring $B = D[x_1,-,x_n]$, where D is a division ring and $x_i$ are commuting variables over D . Theorem-(1.6) ( or theorem -0 ) states the following.

Theorem (1.6) ( Theorem -0 ) :- Let B be the Polynomial ring $B = D[x_1,-,x_n]$, where D is a division ring and $x_i$ are commuting variables over D. Let M be a finitely generated B-module. Then the following hold true ;

(a) If M is a torsion B-module, then there exists a change of variables $t_i$ of $x_i$ such that B can also be written as a polynomial ring $B = D[t_1,-,t_n]$, in the new variables $t_i$ and such that M is a finitely generated $B_{n-1}$-module, where $B_{n-1}$ is the polynomial ring $B_{n-1} = D[t_1,\ldots,t_{n-1}]$ in the fewer new 'n-1' variables.

(b) In case M is further a torsion $B_{n-1}$-module ($B_{n-1}$ is as in (a) above) then there exists a further change of variables $z_i$ of $t_i$; $1 \le i \le n-1$, such that $B= D[z_1,\ldots z_{n-1}]$ and M is a finitely generated $B_{n-2}$-module where $B_{n-2}$ is the polynomial ring $B_{n-2} = D[z_1,-,z_{n-2}]$. Moreover in this case we can write B also as the polynomial ring $B = D[z_1,-,z_n]$, where $z_n$ can be chosen such that $z_n=t_n$ .

Next, we state our theorem-I which makes the statement of theorem-0 that we have stated above more precise by connecting it with the Krull dimension of the module M. For this however we first introduce a definition called by us as the krull dimension condition .

**Defintion (1.7) (Krull Dimension Condition):-** Let B be a noetherian ring and let A be a a noetherian subring of B. If M is a finitely generated

B-module , we say that M has the Krull dimension condition relative to A if ,$|yA|_A \leq |yB|_B$ , for all y in M.

The next lemma shows that the krull dimension condition is satisfied by modules over noetherian polynomial rings relative to certain noetherian subrings .

**Lemma (1.8) :- (a) Let A be a noetherian ring and l**et B=A[x], be a polynomial ring over A in a commuting variable x . If M is a finitely generated B-module then M satisfies the Krull dimension condition relative to A,that is , for any $y \in M$;

$|yA|_A \leq |yB|_B$ .

**(b)** Let A be a Noetherian ring with $|A|=a$ . Let B= A[$x_1$,-,$x_n$], be a polynomial ring over A in n- commuting variables , $x_i$ . Let B'$_m$ =D[$x_1$,-,$x_m$], be the polynomial ring in 'm' variables; where m is an integer such that $0 \leq m \leq n$, with B'$_0$ =D, and B'$_n$ =B. Let M be a finitely generated B-module . Then the following hold true ;

(i) $0 \leq |M| \leq a+n$.

(ii) M satisfies the Krull dimension condition relative to the subring B'$_m$ . That is for any $y \in M$; $|yB'_m|_{B'm} \leq |yB|_B$ . Also we have that $|yB'_m|_{B'm} \leq |yB'_{m+t}|_{B'm+t}$ , for all t>0.

We now introduce another definition called as the strong krull dimension condition ;

**Defintion(1.9) (Strong Krull Dimension Condition**):- Let B be a noetherian ring and let A be a noetherian subring of B. If M is a finitely generated B-module , we say that M has the strong Krull dimension condition relative to A if ,$|yA|_A = |yB|_B$ ,for all y in M.

**Remark**:- Clearly the above definition means that if a module M over a noetherian ring B satisfies the strong krull dimension condition relative to a noetherian subring A of B , then it must have the following properties ;

(i) $|yB|_B \leq |A|$ ,for all y in M and hence,$|M|_B \leq |A|$.

(ii) Moreover it also means that if a module M over a noetherian ring B satisfies the strong krull dimension condition relative to a noetherian subring A of B then M satisfies the krull dimension condition relative to the ring A.

We now state the following lemma (1.10) below ;

**Lemma(1.10) :-** Let B be a noetherian ring and let A be a noetherian subring of B with $|A|=a$ . If M is a finitely generated B-module with $|M|_B =b$ that is also a finitely generated A-module then the following hold true ;

(i) $|M|_B \leq |M|_A \leq |A|$ . Thus $b \leq a$ . In fact, $|yB|_B \leq |yB|_A \leq |A|$ , for all y in M .

(ii) Moreover M satisfies the krull dimension condition relaive to the subring A if and only if $|yB|_B=|yB|_A$ , for all y in M (hence in this case $|M|_B=|M|_A$ ) .

Remark:- We mention that Part (ii) of the lemma (1.10) stated above does not make the claim , namely, if a module M over a noetherian ring B has the krull dimension condition relative to a noetherian subring A , then M has the strong krull dimension condition relative to A . That would mean then that $|yB|_B=|yA|_A$ , for all y in M in (ii) of the above lemma .

The above lemma can be applied to noetherian polynomial rings to yield the following proposition **;**

**Proposition (1.11)**:- (i) Let A be a noetherian ring and let B=A[x] be the polynomial ring over A in the commuting variable x. If M is a finitely generated B-module that is also a finitely generated A-module then we get that $|M|_B=|M|_A$ . In fact $|yB|_B=|yB|_A$ , for all y in M .

(ii) Let B be the Polynomial ring B = A[$x_1$,-,$x_n$], where A is a noetherian ring and $x_i$ are commuting variables over A. Let M be a finitely generated B-module. Let B'$_m$ = A[$x_1$,-,$x_m$], be the polynomial subring of B . If M is a finitely generated B-module that is also a finitely generated B'$_m$-module then $|M|_B = |M|_{B'm}$ . In fact $|yB|_B=|yB|_{B'm}$ , for all y in M .

Before we proceed further we state the following useful Lemma (1.12) below ;

**Lemma (1.12):-** Let R be a semiprime , right noetherian ring that is right krull homogenous . Let M be a finitely generated module over R. Then $|M| <|R|$ iff M is a torsion R-module.

Lemma (1.12) above allows us to state Theorem (1.13) below which says when can a module M over a semiprime , noetherian ring B be torsion ( or non-torsion ) over a noetherian semiprime subring A of B .

**Theorem (1.13)** :- Let B be a noetherian, semiprime ring and let A be a semiprime, noetherian subring of B with $|A|=a$ . Suppose further A is a krull homogenous ring . Let M be a finitely generated B-module with $|M|=b$ and so that M is also a finitely generated A-module. Then the following hold true ;

(i) If M has the krull dimension condition relaive to the subring A, then b<a , if and only if M is a finitely generated , torsion A-module . ( In fact yB is a finitely generated , torsion A-module for all y in M ) .

(ii ) Again if M has the krull dimension condition then b $\geq$ a ( hence b=a ) , if and only if M is a non-torsion A-module. ( It is not necessarily true that yB is a non- torsion A-module , for all y in M) .

The above results especially lemma (1.10) and theorem (1.13 ) can be applied to the polynomial ring B = D$[t_1,-,t_n]$, over a division ring D in 'n' commuting variables $t_i$, to yield theorem-(1.14) below ;

**Theorem (1.14)** :- Let A be a semiprime, noetherian , krull homogenous ring with |A|=a . Let B be the Polynomial ring B = A$[t_1,-,t_n]$ , over A in n - commuting variables $t_i$ . Let $B_m$ denote the subring of B namely , $B_m$= A$[t_1,-,t_m]$, (o$\leq$m$\leq$n), with $B_0$=A and $B_n$=B. Let M be a finitely generated B-module with |M|=b so that M is also a finitely generated $B_m$ - module. Then the following hold true ;

(i) $|M|_B = |M|_{Bm}$ =b and 0 $\leq$ b $\leq$ a+m =|Bm| .

(ii) Also b< |Bm|=a+m , if and only if M is a finitely generated , torsion $B_m$ - module .

(iii) Moreover b $\geq$ a+m ( hence b=a+m ) , if and only if M is a non-torsion $B_m$ -module .

Theorem (1.14) now allows us to apply theoarem-0 to the polynomial ring B = D$[x_1,-,x_n]$, over a division ring D in 'n' commuting variables $x_i$ , to get theorem(1.15) ( named as Theorem-I ) below that characterises the krull dimension of a finitely generated module over the polynomial ring $B = D[x_1,-,x_n]$ .

**Theorem (1.15) ( Theorem -I)** :- Let $B = D[x_1,-,x_n]$, be a polynomial ring in the commuting variables $x_i$ over a division ring D. Let M be a finitely generated B-module. Let $B'_m = D[x_1,-,x_m]$, be the polynomial ring in 'm' variables; where m is an integer such that $0 \leq m \leq n$, with $B'_0 = D$, and $B'_n = B$. Then the following conditions on M are equivalent ;

(i) Krull dimension (M) = m; $0 \leq m \leq n$ .

(ii) There exists a change of variables $t_i$ of $x_i$ such that B can also be expressed as a Polynomial ring $B = A[t_1,-,t_n]$, in the new variables $t_i$ ( giving an automorphism f of B such that $f(x_i)=t_i$ ) and M is a finitely generated non-torsion $B_m$ -module where $B_m$ is the polynomial ring $B_m = D[t_1,-,t_m]$, in the fewer, m variables , $t_1,-,t_m$ ( clearly $f(B'_m)=B_m$ ) . (Moreover , in this case then for any positive integer k , with k>m and $0 \leq k \leq n$ , M is a finitely generated torsion $B_k$ module where $B_k$ is the polynomial ring $B_k = D[t_1,-,t_k]$, in the fewer, k variables , $t_1,-,t_k$ ) .

Now a straight-forward version of theorem (1.15) stated above yields our main theorem given below ;

Theorem(1.16) (Main Theorem) :- Let $B = D[x_1,-,x_n]$, be a polynomial ring in the commuting variables $x_i$ over a division ring D. Let M be a finitely generated B-module. Let $B'_m = D[x_1,-,x_m]$, be the polynomial ring in 'm' variables; where m is an integer such that $0 \leq m \leq n$, with $B'_0 = D$, and $B'_n = B$. Then the following conditions on M are equivalent ;

(i) Krull dimension (M) = m; $0 \leq m \leq n$ .

(ii) M is a non -torsion $B'_m$ -module (not necessarily finitely generated as a $B'_m$ -module) and such that for any k>m , M is a torsion $B'_k$ - module .

**Section(2) (Generalisations)** :- We mention briefly that the above theorem (1.16) can be rephrased more generally for the polynomial ring $B = A[x_1,-,x_n]$, in n-commuting variables $x_i$, over an artinian ring A which makes it possible to further generalise this theorem for the ring $B = A[x_1,-,x_n]$ , in case A is any noetherian ring that has an artinian quotient ring . We state these generalsations below . However before we announce these generalisations recall that in section (1) above , for a noetherian ring R, we denoted by C(0) the set of regular elements of R. Then for a module M over R we gave the definition of a torsion and a torsion-free element as well as the definition of a torsion and a torsion-free sub-module of the module M. Then we stated our key theorems in terms of these concepts . In this section we will first introduce similar definitions as mentioned above and in terms of which we then state our generalisations .

**Definition (2.1):-** Given a noetherian ring R with nilradical N , denote by C(N) the set of elements of R that are regular modulo N. Let M be a module over R . An element , y of M is said to be an N- torsion element of M if yd = 0 for some element d of C(N). An element y of M is said to be an N- torsion-free element of M if $yd \neq 0$ for all elements d of C(N) . Clearly y is a torsion-free element implies that y is a nonzero element of M. If N=(0) , then R is a semiprime ring   and  in  this  case  C(0)   denotes   the   set   of
regular  elements  of  R .

**Definition (2.2):-** Let   R   be   a   noetherian   ring   with
nilradical  N . We   say   a   module   M   over   the
noetherian ring  R  is  an  N-torsion  module  if  every
element   of  M  is  an  N- torsion  element . We  say  the
module   M  over  the  noetherian ring  R  is  an  N- torsion -
free  module  if  every nonzero  element  of  M  is  an  N-
torsion-free  element.

**Remark :-** Let  R  be  a  right   noetherian   ring   with
nilradical  N  and  let  M  be  a   finitely   generated
module over  R . Let   T(M)   be   the   set   of  N- torsion
elements of M . Then  T(M)  is  a  non-empty  subset  of  M
and   T(M)   is   usually   not   a   submodule   of   M .
However  T(M)  is  a  submodule  of  M  if  C(N)  is  a
right  ore  set . It  is  clear  if  N=(0) , then  R   is  a
semiprime   right   Noetherian   ring   and   then   by
Goldie's  Theorem (see  [2]) C(0)  is   a   right  ore  subset
of   R , where   C(0)   is   now  the   set   of   regular
elements  of   R . In   general , from [4] , it  follows  that
C(N) is   a   right  ore  set  in  R   if  and  only  if
C(N)=C(0) , where   C(0)   as  usual  is  the  set  of regular
elements  of  R . Moreover  C(N)=C(0)  if  and  only  if   R
has  an  artinian   quotient  ring .

We  now  mention  first  that  results  similar to  the  results  of  section (1)
can  be  applied  to  the  Polynomial ring   B = A[$x_1$,-,$x_n$],  where   A  is an Artinian ring and $x_i$ are commuting indeterminates over A. These results then culminate in the following theorem for the Polynomial ring B= $A[x_1,-,x_n]$ , where A is an artinian ring .

**Theorem(2.3)** :- Let A be an artinian ring with nilradical N. Let s=C(N) denote the set of elements of A that are regular modulo N . Let B = $A[x_1,-,x_n]$, be a polynomial ring over A in 'n' commuting variables $x_i$. Let $B'_m$ =$A[x_1,-,x_m]$, be the polynomial ring in 'm' variables; where m is an integer such that $0 \leq m \leq n$, with $B'_0$ =A, and $B'_n$ =B . Denote by $N(B'_m)$ the nilradical of the ring $B'_m$ . Let M be a finitely generated B-module . Then the following hold true ;

(a) M is a s -torsion-free A module (M need not be a finitely generated A-module ). Also If $|M| = m$ , for some m, then $0 \leq m \leq n$ .

(b) Moreover the following are equivalent ;

(i) M is a finitely generated B-module with $|M|=m$ , for some m , $0 \leq m \leq n$ .

(ii) M is a $N(B'_m)$-non -torsion $B'_m$ -module (not necessarily finitely generated as a $B'_m$ -module) and such that for any k>m , M is a $N(B'_k)$- torsion $B'_k$ - module .

A further generalisation of Theorem(2.3) is the statement of theorem (2.4) below ;

**Theorem (2.4)** :- Let A be a noetherian ring that has an artinian quotient ring . Let $|A|=a$ and let N be the nilradical of A such that A/N is a krull homogenous ring . Let   B = A[$x_1$,-,$x_n$],  be a polynomial ring over A in 'n' commuting variables $x_i$ . Let  B'$_m$  =A[$x_1$,-,$x_m$], be the polynomial subring oof  B  in 'm' variables; where m is an integer such that $0 \leq m \leq n$, with  B'$_0$  =A, and B'$_n$ =B .  Denote by N(B'$_m$ ) the  nilradical  of  the  ring   B'$_m$  and  let   s= C(N) , denote the  set  of  elements  of  A  that  are  regular  modulo N . Let M  be a finitely  generated   B-module with |M|=b .  Then the following hold  true ;

 (a) If  M  is an  s-torsion-free  B- module with |M|=b, then a $\leq$ b $\leq$ a+n , hence  |M|=b=a+m , for some m with  $0 \leq m \leq n$ .

(b) Moreover  the  following  are  equivalent ;

(i) M  is  a  finitely  generated ,  B - module   with    |M|=a+m .

(ii) M  is  an N(B'$_m$)-non -torsion  B'$_m$- module (not  necessarily finitely  generated  as  a  B'$_m$ -module) such  that   for  any positive  integer k, with   k>m ,   M  is  an  N(B'$_k$) -torsion module .

**Remark**:- (a) There exists  a   simple  noetherian  domain  with |A|=1 with  a  simple  module  M ( clearly  M is  a  faithful  A  module ). Let  s denote  the  set  of  nonzero elements    of  A .   Let  B  denote  the polynomial  ring   B = A[$x_1$,-,$x_n$],   over  A  in n-commuting variables $x_i$ . Let M'=M[$x_1$,$x_2$,...,$x_n$] . Then  M'  is a finitely generated  B-module  such that  M'  is  a s-torsion- B- module . Thus    the  stated hypothesis  of theorem (1.16)   does  not  hold true  for  the  module  M'.

(b) If   A  is  a  commutative  noetherian  ring  and   B = A[$x_1$,-,$x_n$], is the  polynomial  ring  over  A  in n-commuting variables $x_i$ , then  for any  finitely generated ,  B- module  M'  such  that   M'  is  a   primary faithful A- module the hypothesis and hence the statement of theorem(1.16) always holds true .

## References:-